\documentclass[12pt,a4paper]{amsart}
\usepackage{a4wide}
\usepackage[english]{babel}
\usepackage[T2A]{fontenc}
\usepackage[utf8]{inputenc} 
\usepackage{amsfonts}
\usepackage{amssymb, amsthm, amscd}
\usepackage{amsmath}
\usepackage{mathtools}
\usepackage{needspace}
\usepackage{etoolbox}
\usepackage{lipsum}
\usepackage{comment}
\usepackage{cmap}
\usepackage[pdftex]{graphicx}
\usepackage[unicode]{hyperref}
\usepackage[matrix,arrow,curve]{xy}
\usepackage[usenames,dvipsnames]{xcolor}
\usepackage{colortbl}
\usepackage{textcomp}
\usepackage{cite}
\usepackage{euscript}

\DeclareMathOperator{\pr}{pr}

\newcommand{\Z}{\mathbb{Z}}

\newcommand{\eps}{\varepsilon}
\newcommand{\ph}{\varphi}

\newcommand{\sub}{\subseteq}

\newcommand{\Ker}{\mathop{\mathrm{Ker}}\nolimits}

\renewcommand{\ge}{\geqslant}
\renewcommand{\le}{\leqslant}
\newcommand{\sm}{\setminus}
\newcommand{\map}[3]{#1\colon #2\to #3}

\renewcommand{\>}{\rangle}

\newcommand{\SR}{\mathrm{SR}}
\newcommand{\ASR}{\mathrm{ASR}}

\theoremstyle{plain}
\newtheorem{thm}{Theorem}
\newtheorem{lem}{Lemma}

\theoremstyle{definition}

\newtheorem{cor}{Corollary}

\theoremstyle{remark}
\newtheorem*{rem}{Remark}

\AtBeginEnvironment{thm}{\begin{samepage}}
	\AtEndEnvironment{thm}{\end{samepage}}
\AtBeginEnvironment{lem}{\begin{samepage}}
	\AtEndEnvironment{lem}{\end{samepage}}
\AtBeginEnvironment{st}{\begin{samepage}}
	\AtEndEnvironment{st}{\end{samepage}}
\AtBeginEnvironment{crit}{\begin{samepage}}
	\AtEndEnvironment{crit}{\end{samepage}}
\AtBeginEnvironment{ax}{\begin{samepage}}
	\AtEndEnvironment{ax}{\end{samepage}}
\AtBeginEnvironment{defn}{\begin{samepage}}
	\AtEndEnvironment{defn}{\end{samepage}}
\AtBeginEnvironment{cor}{\begin{samepage}}
	\AtEndEnvironment{cor}{\end{samepage}}
\AtBeginEnvironment{note}{\begin{samepage}}
	\AtEndEnvironment{note}{\end{samepage}}
\AtBeginEnvironment{prop}{\begin{samepage}}
	\AtEndEnvironment{prop}{\end{samepage}}

\newcommand{\SL}{\mathop{\mathrm{SL}}\nolimits}

\newcommand{\tc}{\text{,}}
\newcommand{\tp}{\text{.}}
\renewcommand{\tilde}{\widetilde}

\newcommand{\vpi}{\varpi}

\DeclareMathOperator{\Max}{Max}

\makeatletter
\def\@settitle{\begin{center}%
		\baselineskip14\p@\relax
		\bfseries
		\@title
	\end{center}%
}

\def\@evenhead{\hfil\sc Pavel Gvozdevsky\hfil}
\def\@oddhead{\hfil\sc Improved $K_1$-stability\hfil}
\makeatother

\title{Improved $K_1$-stability for the embedding $D_5$ into $E_6$}

\author{Pavel Gvozdevsky}
\date{\today}
\address{Chebyshev Laboratory, St. Petersburg State University, 14th Line V.O., 29B, Saint Petersburg 199178 Russia}
\thanks{Research is supported by "Native towns", a social investment program of PJSC "Gazprom Neft".}

\begin{document}

\maketitle

\begin{abstract}
	This paper is dedicated to the surjective stability of the $K_1$-functor for Chevalley groups for the embedding $D_5$ into $E_6$. This case was already studed by Plotkin. In the present paper, we improve his result by showing that surjective stability holds under a weaker assumption on a ring. Another result of the present paper shows how the $K_1$-stability can help to study overgroups of subsystem subgroups.
\end{abstract}

\section{Introduction}
The problem of $K_1$-stability for semisimple algebraic groups was first discussed by Hyman Bass in his book \cite{BassBook}. He suggested that stability hypotheses can be formulated in terms of Jacobson dimension of a ring and the relative semisimple rank of a group. 

In the framework of Chevalley groups, this problem was first studied by Stein (see \cite{SteinStability}, \cite{SteinMats}), who developed an approach that allows to consider both classical and exceptional cases in the same style. This problem consists of two independent parts, surjective and injective stability.

In 1999 Vavilov suggested that injective stability of the $K_1$-functor for all regular embeddings of root systems holds under natural conditions on the stable rank of a ring, depending on embedding. This result is based on the Suslin--Tulenbaev proof (see \cite{SusTul}) of the Dennis--Vaserstein decomposition (see \cite{VasersteinUniStab}). For classical groups, this idea was implemented by A. Bak, V. Petrov, and G Tang (see \cite{BakTang} and \cite{BakPetrTang}).

The surjective stability question was studied by Stein in \cite{SteinStability}, by Vaserstein in \cite{VasersteinStability} and \cite{VasersteinUniStab}, and by Plotkin in \cite{PlotProc},\cite{PlotkinStability} and \cite{PlotkinE7}. In those results there are several types of conditions on a ring. Depending on the case it may be a condition on the stable rank, on the absolute stable rank, on the dimension of maximal spectrum, the condition $V_n$, introduced by Vaserstein in \cite{VasersteinUniStab} for orthogonal and unitary groups, or the similar condition $V_\mu(\Phi,R)$, introduced by Plotkin in \cite{PlotkinE7} for certain representation of exceptional groups. In particular, according to \cite{PlotkinStability}, the surjective stability of the $K_1$-functor for the embedding $D_5$ into $E_6$ holds when the ring satisfies the condition $\ASR_4$,below we recall the definition. In the present paper, we show that one can replace $\ASR_4$ here by a weaker condition $\ASR_5$. 

In the last section, I make an observation how the $K_1$-stability can help to study overgroups of the subsystem subgroup $E(A_1+D_6,R)$ in $G(E_7,R)$.  

\section{Preliminaries and notation}

\subsection{Chevalley groups and $K_1$-functor}
Let $\Phi$ be a reduced irreducible root system of rank $l$, let $G(\Phi,-)$, be a simply connected Chevalley--Demazure group scheme over $\Z$ of type $\Phi$ (see \cite{ChevalleySemiSimp}), and let $T(\Phi,-)$, be a split maximal torus in it. If $R$ is a commutative ring with unit, the group $G(\Phi,R)$ is called the simply connected Chevalley group of type $\Phi$ over $R$.

For a subset $X\sub R$ we denote by $\<X\>$ the ideal generated by $X$, similarly for a subset $X$ of a group we denote by $\<X\>$ the subgroup generated by $X$.

 To each root $\alpha\in\Phi$ there correspond root unipotent elements $x_\alpha(\xi)$, $\xi\in R$,  elementary with respect to $T$. The group generated by all these elements 
 $$
 E(\Phi,R)=\left\<x_\alpha(\xi)\colon\alpha\in\Phi\tc\,\xi\in R\right\>\le G(\Phi,R)\tp
 $$
 is called the elementary subgroup of $G(\Phi,R)$.
 
It is well known that if $\Phi$ is an irreducible root system of rank $l\ge 2$, then $E(\Phi,R)$ is always normal in $G(\Phi,R)$, see, for example, \cite{Taddei}. 

Thus, the $K_1$-functor of Chevalley group $G(\Phi,R)$ is naturally defined as the quotient group
$K_1(\Phi,R)=G(\Phi,R)/E(\Phi,R)$, see, for example, \cite{SteinStability}.

Any inclusion of root systems $\Delta\sub\Phi$ induces the homomorphisms $G(\Delta,R)\to G(\Phi,R)$,  $E(\Delta,R)\to E(\Phi,R)$ taking roots to roots, and the homomorphism of the corresponding $K_1$-functors $K_1(\Delta,R)\to K_1(\Phi,R)$. The surjective stability problem is to find sufficient conditions on the ring $R$, depending on $\Delta$ and $\Phi$, for the homomorphism on $K_1$ to be surjective. In other words, to find sufficient conditions for the existence of the decomposition $G(\Phi,R)=E(\Phi,R)G(\Delta,R)$.

\subsection{ASR-condition}
Recall that a commutative ring $R$ satisfies the absolute stable rank condition $\ASR_n$ if for any row $(r_1,\ldots,r_n)$ with coordinates in $R$, there exist elements $t_1$,$\ldots$,$t_{n-1}\in R$ such that every maximal ideal of $R$ containing the ideal $\<r_1+t_1r_n,\ldots,r_{n-1}+t_{n-1}r_n\>$ contains already the ideal $\<r_1,\ldots,r_n\>$. This notion was introduced in \cite{EstesOhm} and used in \cite{SteinStability}, \cite{SteinMats} and then in \cite{PlotProc},\cite{PlotkinStability}, and \cite{PlotkinE7} to study stability problems.

If we assume that a row $(r_1,\ldots,r_n)$ is unimodular, then the absolute stable rank condition boils down to the usual stable rank condition $\SR_n$ (see \cite{Bass64},\cite{VasersteinStability}).

Absolute stable rank satisfies the usual properties, namely
for every ideal $I\unlhd R$ condition $\ASR_n$ for $R$ implies $\ASR_n$ for the quotient $R/I$, and if $n\ge m$, then $\ASR_m$ implies $\ASR_n$ . Finally, it is well known that if the dimension of the
maximal spectrum $\dim\Max(R)$ is $n-2$, then both conditions $\ASR_n$
and $\SR_n$ are satisfied (see \cite{EstesOhm},\cite{MagWandarKalVas},\cite{SteinStability}).

\subsection{Weight diagrams and basic representation}

Let us fix an order on $\Phi$, and let $\Phi^+$, $\Phi^-$ and $\Pi=\{\alpha_1\tc\ldots\tc\alpha_l\}$ be the sets of positive, negative, and fundamental roots, respectively. Our numbering of the fundamental roots follows that of \cite{Bourbaki4-6}. By $\vpi_1$,$\ldots$,$\vpi_l$ one denotes the corresponding fundamental weights. Let $W=W(\Phi)$ be the Weyl group of the root system $\Phi$.

Recall that an irreducible representation $\pi$ of the complex semisimple Lie algebra $L$ is called {\it basic} (see \cite{Matsumoto69}) if the Weyl group $W=W(\Phi)$ acts transitively on the set $\Lambda(\pi)$ of nonzero weights of the representation $\pi$. 

In this paper we can restrict ourselves to the basic representations without zero weight. In this case, $\Lambda(\pi)$ is the set of all weights of $\pi$, and the weights $\Lambda(\pi)$ form the one Weyl orbit. Such representations are called {\it microweight} or {\it minuscule} representations, and the list of these representations is classically known (see \cite{Bourbaki7-9})

Let $L$ be a simple Lie algebra of type $\Phi$. To each complex representation $\pi$ of the algebra $L$ there corresponds a representation $\pi$ of the Chevalley group $G=G(\Phi,R)$ on the free $R$-module $V=V_R=V_Z\otimes_{\Z}R$ (see \cite{Matsumoto69},\cite{Steinberg85}). If $\pi$ is faithful we can identify $G$ with its image $\pi(G)=G_{\pi}(\Phi,R)$ under this representation and omit the symbol $\pi$ in the action of $G$ on $V$. Thus, for an $x\in G$ and $v\in V$ we write $xv$ for $\pi(x)v$. If we want to specify that the group $G=G(\Phi,R)$ is considered in the basic representation $\pi$ with the highest weight $\mu$, the notation $(G(\Phi,R),\mu)$ is used. In the sequel, $\mu$ always stands for the highest weight of a representation.

Decompose the module V into the direct sum of its weight submodules
$$
V=\bigoplus_{\lambda\in\Lambda(\pi)} V^{\lambda}\tp
$$
It follows from the definition of microweight representations that all $V^{\lambda}$, $\lambda\in\Lambda(\pi)$, are one-dimensional. Matsumoto, \cite{Matsumoto69} Lemma 2.3, has shown that there is a special base of weight vectors $e^\lambda\in V^\lambda$, $\lambda\in\Lambda(\pi)$, such that the action of root unipotents $x_{\alpha}(\xi)$, $\alpha\in\Phi$, $\xi\in R$, is described by the a following simple formulas:
\begin{equation}
\begin{split}
&\text{i.  if }\lambda\in\lambda(\pi)\tc\, \lambda+\alpha\notin\Lambda(\pi)\text{, then }x_{\alpha}(\xi)e^\lambda=e^{\lambda}\tc\\
&\text{ii.  if }\lambda\in\lambda(\pi)\tc\, \lambda+\alpha\in\Lambda(\pi)\text{, then }x_{\alpha}(\xi)e^\lambda=e^{\lambda}\pm\xi e^{\lambda+\alpha}\tp
\end{split}\notag
\end{equation}

Now we may conceive any element $g\in G$ as a matrix $g=(g_{\lambda,\nu})$, where $\lambda$ and $\nu$ range over all the weights of the representation $\pi$. We will denote the $\nu$-th column and the $\lambda$-th row of this matrix by $g_{*,\nu}$ and $g_{\lambda,*}$, respectively.

Weight diagrams serve as a great visual aid for calculations in Chevalley groups (see \cite{atlas},\cite{SteinStability},\cite{VavThirdlook} for the detailed description of these diagrams and references). Let us recall here the corresponding definitions.

It is well known that a choice of a fundamental system $\Pi$ defines a partial order of the weight lattice $P(\Phi)$ as follows: $\lambda\ge\nu$ if and only if $\lambda-\nu$ is a linear combination of the fundamental roots with nonnegative integer coefficients. Let us associate with a representation $\pi$ a graph, which is in fact the Hasse diagram for the set $\Lambda(\pi)$ of its weights with respect to the above order.

We construct a labeled graph as follows. Its vertices correspond to the weights $\lambda\in\Lambda(\pi)$ of the representation $\pi$, and the vertex corresponding to $\lambda$ is actually labeled with $\lambda$ (usually these labels are omitted).

We read the diagram from right to left and from bottom to top, which means that a larger weight tends to stand to the left of and higher than a smaller one. The leftmost vertex corresponds to the highest weight $\mu$ of a representation.

The vertices corresponding to the weights $\lambda$,$\nu\in\Lambda(\pi)$ are joined by a bond labeled with $\alpha_i$ (or simply i) if and only if their difference $\lambda-\nu=\alpha_i\in \Pi$ is a fundamental root. We draw the diagrams in such a way that the labels on the opposite sides of a parallelogram are equal (as a rule, at least one of them is omitted). Thus, all paths of minimal length connecting two vertices have the same sum of labels. This means that if there is a root, corresponding to a pair of vertices, it can be determined by any path of minimal length between these vertices.

Now one may conceive a vector $v\in V$ as such a weight diagram which has an element $v^{\lambda}\in R$ attached to every node. The above-mentioned Matsumoto lemma gives a very simple rule describing what happens with such a vector $v$ under the action of $x_\alpha(\xi)$. For a minuscule $\pi$ and a fundamental root $\alpha_i$, $x_{\alpha_i}(\xi)$ adds or subtracts (always adds for a clever choice of the weight base) $\xi v^\nu$ to $v^\lambda$ along each edge labeled with $i$. For other roots one merely has to trace all paths in the diagram which have the same labels at their edges as the root $\alpha$ in its linear expansion with respect to the fundamental roots. The order of the labels on such a path together with the structure constants of the Lie algebra determines the sign with which $x_\alpha(\xi)$ acts on $v^\lambda$ (see \cite{VavSigns}).

In this paper we use only microweight representations; that is, we do not care about weight diagrams with zero weights. The details of how to construct and operate with weight diagrams in the presence of zero weight can be found in \cite{atlas}. Weight diagrams as a tool for proving $K_1$ and $K_2$-functors stability were introduced by Stein in \cite{SteinStability}.

\subsection{Chevalley--Matsumoto Decomposition}

Stability of the $K_1$-functor is closely related to the fact known as the Chevalley--Matsumoto decomposition theorem (see \cite{ChevalleySemiSimp},\cite{Matsumoto69},\cite{SteinStability}). Let us formulate the special case of this theorem used for the stability problem. Consider a basic representation $\pi$ of the group $G(\Phi,R)$ with the highest weight $\mu$. Denote by $\alpha_k\in\Pi$ the fundamental root such that $\mu-\alpha_k$ is a weight, and by $\Delta$ the subsystem in $\Phi$ generated by all fundamental roots except $\alpha_k$. Further, let $\Sigma=\Phi\sm\Delta$, $\Sigma^+=\Phi^+\cap \Sigma$, $\Sigma^-=\Phi^-\cap \Sigma$, and
\begin{gather*}
U(\Sigma, R)=\<x_\alpha(\xi)\tc\, \alpha\in \Sigma^+\tc\, \xi\in R\>\tc\\
V(\Sigma, R)=\<x_\alpha(\xi)\tc\, \alpha\in \Sigma^-\tc\, \xi\in R\>\tp
\end{gather*}
Now, take a matrix $g\in (G(\Phi,R),\mu)$ and suppose that $g_{\mu,\mu}$ is an invertible element in $R$. Then the Chevalley--Matsumoto theorem states that $g$ can be expressed in the form
$$
g=vg_1u\tc
$$
where $v\in V(\Sigma, R)$, $u\in U(\Sigma, R)$, and $g_1\in T(\Phi,R)G(\Delta,R)$. 

In particular, $g\in E(\Phi,R)G(\Delta,R)$. Thus, the Chevalley--Matsumoto theorem yields that to prove the surjective $K_1$-stability for the embedding $\Delta\sub\Phi$ it is enough to get a unit of the ring in the left corner of the matrix $g\in (G(\Phi,R),\mu)$ by elementary transformations.

\section{Stability theorem for $D_5\sub E_6$}

\begin{figure}[h]
	\begin{center}
		\includegraphics[scale=0.5]{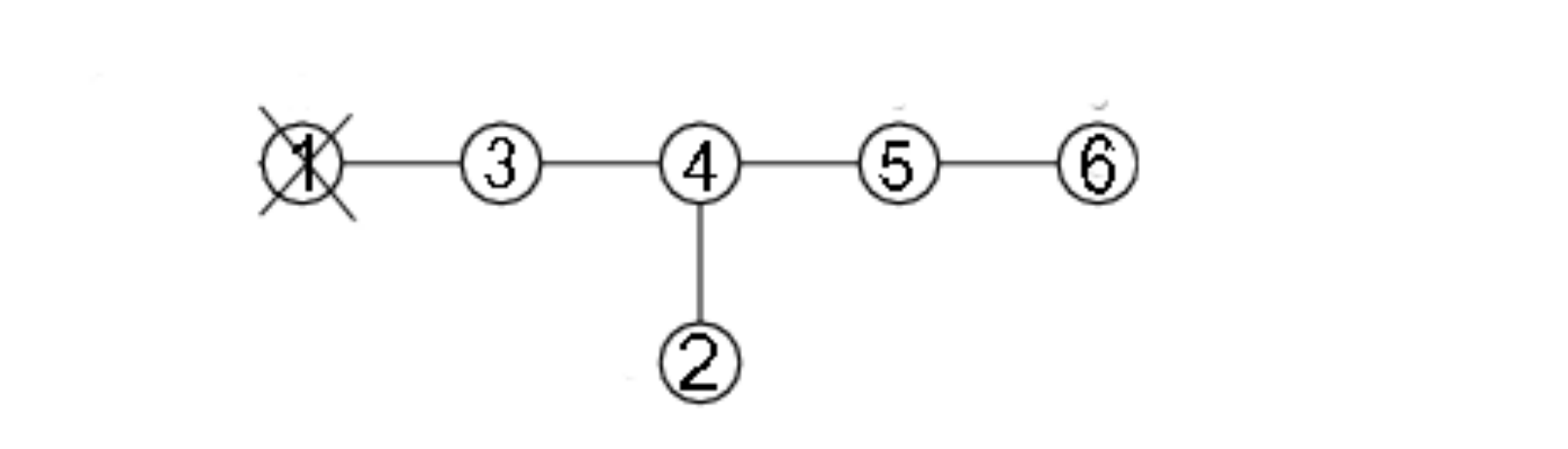}
	\end{center}
	\caption{}
	\label{DynkinE6}
\end{figure}

Consider the inclusion of root systems $D_5\sub E_6$, where the subsystem $D_5$ is generated by all fundamental roots of the system $E_6$ except $\alpha_1$ (see Figure \ref{DynkinE6}).

Our objective is to prove the following theorem.

{\thm\label{D5inE6} Let $R$ be a commutative ring satisfying $\ASR_5$. Then the map
$$
K_1(D_5,R)\to K_1(E_6,R)
$$
is surjective.}

\begin{figure}[h]
	\begin{center}
		\includegraphics[scale=0.4]{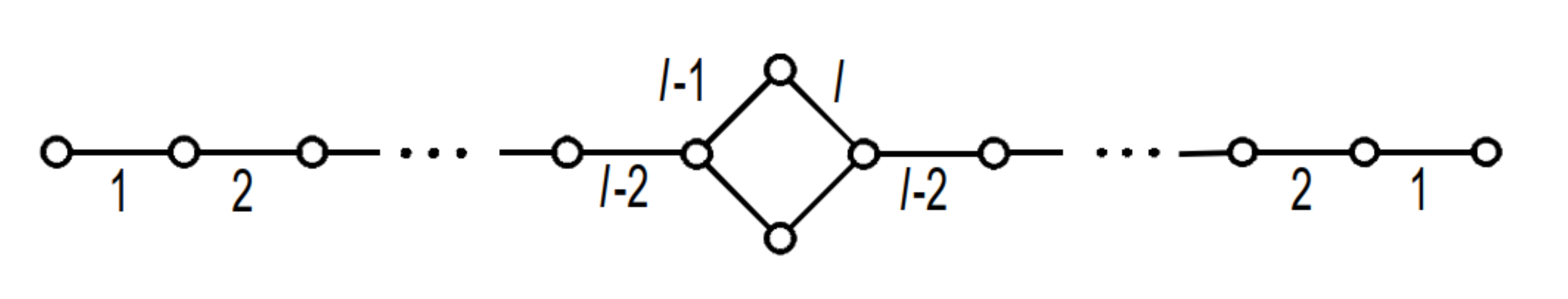}
	\end{center}
	\caption{}
	\label{WeightDiagramDl}
\end{figure}

{\lem\label{SteinDl} Let the ring $R$ satisfy $\ASR_l$, and let $v$ be a unimodular vector in the representation of the group $G(D_l,R)$ with the highest weight $\mu=\vpi_1$, i.e. an unimodular column with coordinates indexed by vertices of the weight diagram at Figure {\rm\ref{WeightDiagramDl}}. Then there exists an element $h\in E(D_l,R)$ such that $(hv)^{\mu}=1$.}

\begin{proof}
	That follows from the proof of $D_l$-case in Theorem 2.1 of \cite{SteinStability}.
\end{proof}

Due to the Chevalley--Matsumoto decomposition, the following lemma implies Theorem \ref{D5inE6}.

\begin{figure}[h]
	\begin{center}
		\includegraphics[scale=0.5]{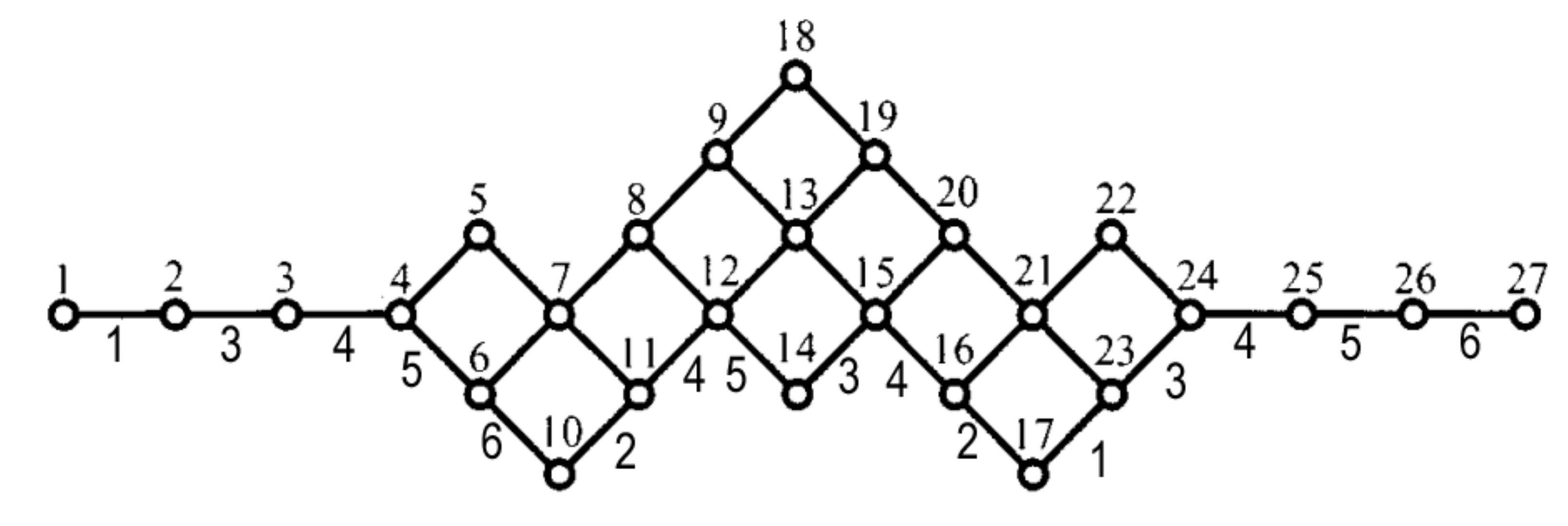}
	\end{center}
	\caption{}
	\label{WeightDiagramE6}
\end{figure}

{\lem Let the ring $R$ satisfy $\ASR_5$, and let $v$ be a unimodular vector in the representation of the group $G(E_6,R)$ with the highest weight $\mu=\vpi_1$ {\rm (}Figure \ref{WeightDiagramE6}{\rm )}. Then there exists an element $h\in E(E_6,R)$ such that $(hv)^{\mu}=1$.}

\begin{proof}
	Using numbering at Figure \ref{WeightDiagramE6} we will write the vector $v$ as a row $(v^1,\ldots,v^{27})$.
\smallskip	

	{\bf Step 1.} Make the row $(v^2,\ldots,v^{27})$ unimodular.	
\smallskip		

	Indeed, let $I=\<v^5,v^7,v^8,v^9,v^{11},\ldots,v^{27}\>\unlhd R$. Applying $\ASR_6$, we can add multiples of $v^1$ to $v^2$,$v^3$,$v^4$,$v^6$,$v^{10}$ as in the definition of $\ASR$-condition. It is easy to see that this additions use only roots from the span of $\alpha_1$,$\alpha_3$,$\ldots$,$\alpha_6$; hence the ideal generated by $v^5$,$v^7$,$v^8$,$v^9$,$v^{11}$,$\ldots$,$v^{27}$ does not change. Now, if a maximal ideal contains new $v^2$,$\ldots$,$v^{27}$, then it in particular contains new $v^2$,$v^3$,$v^4$,$v^6$,$v^{10}$; hence it contains old $v^1$,$v^2$,$v^3$,$v^4$,$v^6$,$v^{10}$, and it also contains the ideal $I$, which means that the initial row was not unimodular. 
\smallskip
		
	{\bf Step 2.} Make the row $(v^1,v^{18}\ldots,v^{27})$ unimodular.
\smallskip
		
	Indeed, let $I=\<v^{18},\ldots,v^{27}\>\unlhd R$. Modulo $I$, the row $(v^2,\ldots,v^{17})$ is unimodular. Hence with the corresponding root elements we can add multiples of $v_2,\ldots,v_{17}$ to $v_1$ such a way that $v_1$ becomes congruent to $1$ modulo $I$. Simultaneously, there will be additions of $v_{18},\ldots,v_{27}$ to $v_2,\ldots,v_{17}$, but modulo $I$ they do not change anything. 
\smallskip
		
	{\bf Step 3.} Make the row $(v^1,v^{18})$ unimodular.
\smallskip
		
	Indeed, let $I\<v_1\>\unlhd R$. Modulo $I$, the row $(v_{18},\ldots,v_{27})$ is unimodular, and the group $E(D_5,R)$ acts on it as in the representation in Lemma \ref{SteinDl}. By the above mentioned properties of the $\ASR$-condition, the quotient $R/I$ satisfies $\ASR_5$. So we just apply Lemma \ref{SteinDl}.
\smallskip
		
	{\bf Step 4.} Make $v_1$ a unit.
\smallskip
		
	Since the row $(v_1,v_{18})$ is unimodular, so is the row $(v_1,\ldots,v_9, v_{18})$. Thus, we can apply Lemma \ref{SteinDl} again but to another copy of $D_5$; namely, the one spanned by all fundamental roots except $\alpha_6$.
\end{proof}

{\rem Actually, as noticed by Nikolai Vavilov, for the first two steps it is enough to require $\ASR_{17}$.}

{\cor\label{K1iso} Let $R$ be a ring satisfying $\ASR_5$ and $\SR_4$, and suppose that its space of maximal ideals has dimension less or equal to $4$. Then all homomorphisms in the following diagram induced by root systems embeddings are isomorphisms.
$$
\xymatrix{
	K_1(D_5,R)\ar[r]\ar[d] & K_1(E_6,R)\ar[d]\\
	K_1(D_6,R)\ar[r] & K_1(E_7,R)\tp
}
$$}

\begin{proof}
	The left arrow is surjective whenever the ring $R$ satisfies $\ASR_6$ (Theorem 2.2 of \cite{SteinStability}), and injective whenever $\SR_5$ holds (Corollary 3.2 of \cite{SteinStability}). The right and upper arrows are injective under condition $\SR_5$ or $\SR_4$ respectively (Theorem 4.1 of \cite{SteinStability}). The right arrow is surjective if the space of maximal ideals has dimension less or equal to $4$ (Theorem 1 of \cite{PlotkinE7}). The upper row is surjective whenever $\ASR_5$ holds (our Theorem \ref{D5inE6}). Finaly, the lower arrow is bijection because the diagram is commutative. 
\end{proof}

{\rem All the above conditions imposed on a ring $R$ in Corollary \ref{K1iso} are satisfied if the ring $R$ is finitely generated, and has Krull dimension less or equal to $3$. Indeed, firstly, the space of maximal ideals has dimension less or equal to the Krull dimension, which on its turn by hypothesis is less or equal to $3$, and that implies $\ASR_5$ (Theorem 2.3 of \cite{EstesOhm}). Secondly, since the ring $R$ is finitely generated, the condition on the Krull dimension also implies $\SR_4$ (Theorem 18.2 of \cite{VasSusSerre}).}

\section{Application to overgroups of $E(A_1+D_6,R)$ in $E(E_7,R)$}

\begin{figure}[h]
	\begin{center}
		\includegraphics[scale=0.5]{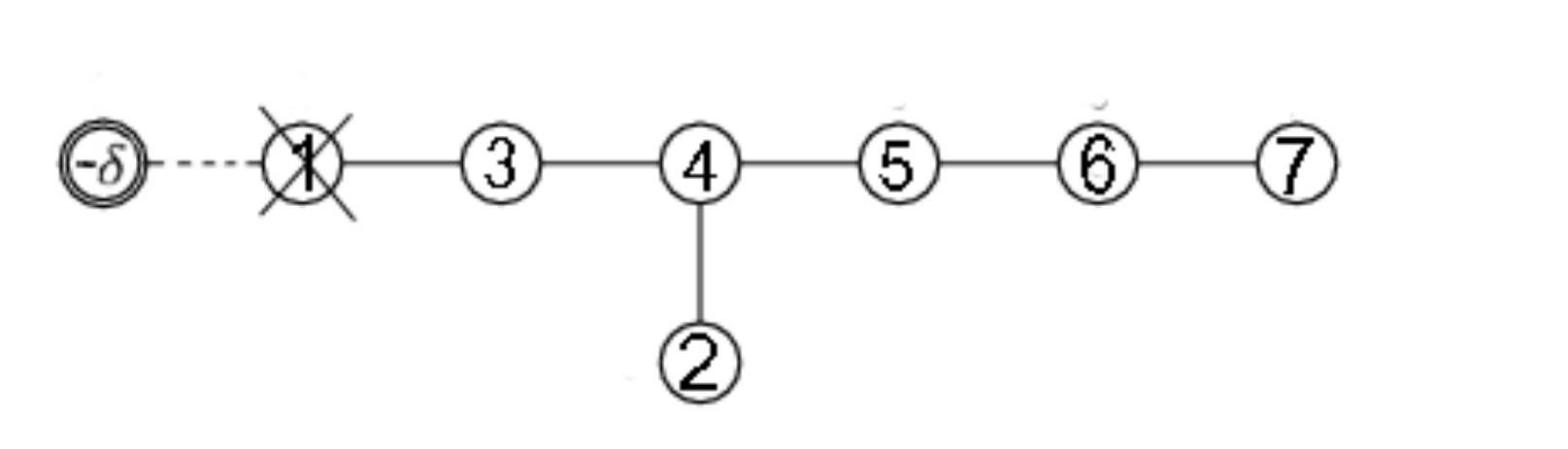}
	\end{center}
	\caption{}
	\label{DynkinDiagramE7}
\end{figure}

In this section, let $\Phi=E_7$ and let $\Delta=A_1+D_6\sub\Phi$ be the subsystem generated by maximal root $\delta$ of $E_7$ and all fundamental roots of $E_7$ except $\alpha_1$, see Figure \ref{DynkinDiagramE7}).

As before, will denote by $G(\Phi,R)$ the simply connected Chevalley group. However,  the corresponding subsystem subgroup $G(\Delta,R)$ is not simply connected. As usual let $E(\Delta,R)$ be its elementary subgroup.

For an ideal $I\unlhd R$ we denote by $\rho_I$ the reduction homomorphism
$$
\map{\rho_I}{G(\Phi,R)}{G(\Phi,R/I)}\tc
$$
induced by projection $R\to R/I$.

The kernel of this homomorphism $G(\Phi,R,I)=\Ker(\rho_I)\le G(\Phi,R)$ is called the {\it principal congruence subgroup}. Set $E(\Phi,I)=\<x_\alpha(\xi)\tc\,\alpha\in\Phi\tc\,\xi\in I\>$. Then, by definition, the {\it relative elementary subgroup} $E(\Phi,R,I)$ of level $I$ is the normal closure of $E(\Phi,I)$ in $E(\Phi,R)$, actually it is automatically normal in $G(\Phi,R)$). 

From the Chevalley commutator formula it follows that $E(\Phi,R,I)$ is generated as a subgroup by the elements $z_{\alpha}(\xi,\zeta)=x_{-\alpha}(-\zeta)x_\alpha(\xi)x_{-\alpha}(\zeta)$, where $\alpha\in\Phi$, $\xi\in I$, $\zeta\in R$. This statement was first formulated be Vaserstein in \cite{VasersteinChevalley} by actually follows from calculations made in papers of Stein \cite{SteinChevalley} and Tits \cite{TitsCongruence}.

We also introduce the following notation:
$$
E(\Phi,\Delta,R,I)=\<E(\Delta,R)\cup E(\Phi,I)\>\le G(\Phi,R)\tp
$$

Our objective is to prove the following theorem.

{\thm\label{EE} Let $R$ be a commutative ring with $2\in R^*$, and let $I\unlhd R$ be an ideal. Then
$$
E(E_7,R,I)\le E(E_7,A_1+D_6,R,I)\tp 
$$}

It can be shown that certain weak form of the sandwich classification theorem holds for overgroups of $E(\Delta,R)$ in $G(\Phi,R)$. My paper with the proof of much more general result will appear in Saint Petersburg Math. Journal. Namely, for any commutative ring $R$ and for any intermediate subgroup $E(\Delta,R)\le H\le G(\Phi,R)$ there exists a unique ideal $I\unlhd R$ such that 
$$
E(\Phi,\Delta,R,I)\le H\le \rho_I^{-1}(G(\Delta,R/I))\tp
$$
 The slice of the subgroup lattice between $E(\Phi,\Delta,R,I)$ and $\rho_I^{-1}(G(\Delta,R/I))$ is what we following Bak call sandwich. The result of Theorem \ref{EE} turns out to be redundant for the proof of this classification, but, probably, this result can help to understand the structure of individual sandwich.   

Theorem \ref{EE} reminds Theorem 3.4 of Alexei Stepanov paper \cite{StepComput}, asserting that $E(\Phi,R,I)\le \<E(\Phi,I),U_P(R)\>$, where $U_P(R)$ is the unipotent radical of a parabolic subgroup $P$. However, the proof of Stepanov's theorem uses basically nothing except the Chevalley commutator formula, and I do not see how Theorem \ref{EE} could be proved in a similar manner. Our proof essentialy relies on the $K_1$-stability. Moreover, condition on a ring in Plotkin's result is not weak enough for our purposes, and this is precisely what motivated me to prove Theorem \ref{D5inE6}.   

{\lem\label{moduloxi^2} Let $R=\Z[{1\over 2}][\xi,\zeta]/(\xi^2)$, and $I=(\xi)\unlhd R$. Then $z_{\alpha_1}(\xi,\zeta)\in E(\Phi,\Delta,R,I)$.}
\begin{proof}
	The set $\{\alpha_1,-\alpha_1\}$ is a subsystem of type $A_1$, consider the corresponding map
	$$
	\map{\ph}{\SL(2,R)}{G(E_7,R)}\tp
	$$
	One can notice that
	\begin{align*}
	z_{\alpha_1}(\xi,\zeta)=\ph\left(
	\begin{pmatrix}
	1 & 0 \\
	-\zeta & 1
	\end{pmatrix}
	\begin{pmatrix}
	1 & \xi \\
	0 & 1
	\end{pmatrix}
	\begin{pmatrix}
	1 & 0 \\
	\zeta & 1
	\end{pmatrix}
	\right)=\ph\left(
	\begin{pmatrix}
	1+\zeta\xi & \xi\\
	-\zeta^2\xi & 1-\zeta\xi
	\end{pmatrix}\right)=\\=\ph\left(
	\begin{pmatrix}
	1 & \xi\\
	0 & 1
	\end{pmatrix}
	\begin{pmatrix}
	1 & 0\\
	-\zeta^2\xi & 1
	\end{pmatrix}
	\begin{pmatrix}
	1+\zeta\xi & 0\\
	0 & 1-\zeta\xi
	\end{pmatrix}\right)=x_{\alpha_1}(\xi)x_{-\alpha_1}(-\zeta^2\xi)h_{\alpha_1}(1+\zeta\xi)\tc
	\end{align*}
	where $h_{\alpha}(\eps)\in T(\Phi,R)$ for $\alpha\in\Phi$, $\eps\in R^*$ are the elements defined by the following formulas
	\begin{gather*}
	w_\alpha(\eps)=x_\alpha(\eps)x_{-\alpha}(-\eps^{-1})x_\alpha(\eps)\tc\\
	h_\alpha(\eps)=w_\alpha(\eps)w_\alpha(-1)\tp
	\end{gather*}
	Moreover, it can be shown that, since $\delta=2\alpha_1+2\alpha_2+3\alpha_3+4\alpha_4+3\alpha_5+2\alpha_6+\alpha_7$, we have
	\begin{align*}
	h_\delta\left(1+{\zeta\xi\over 2}\right)=h_{\alpha_1}\left(1+\zeta\xi\right)h_{\alpha_2}\left(1+\zeta\xi\right)h_{\alpha_3}\left(1+{3\zeta\xi\over 2}\right)h_{\alpha_4}\left(1+2\zeta\xi\right)h_{\alpha_5}\left(1+{3\zeta\xi\over 2}\right)h_{\alpha_6}\left(1+\zeta\xi\right)\\h_{\alpha_7}\left(1+{\zeta\xi\over 2}\right)\tp
	\end{align*}
	Therefore, $h_{\alpha_1}(1+\zeta\xi)\in E(\Delta,R)$, and $z_{\alpha_1}(\xi,\zeta)\in E(\Phi,\Delta,R,I)$.
\end{proof}

{\rem It can be shown that if $R=\Z[\xi,\zeta]/(\xi^2)$, then $h_{\alpha_1}(1+\zeta\xi)\notin E(\Phi,\Delta,R,I)$; hence $z_{\alpha_1}(\xi,\zeta)\notin E(\Phi,\Delta,R,I)$. Therefore, the condition $2\in R^*$ is essential.}

{\lem \label{universal} Let $R=\Z[{1\over 2}][\xi,\zeta]$, and $I=(\xi)\unlhd R$. Then
	\begin{enumerate}
		\item $E(\Phi,R,I^2)\le E(\Phi,I)$.
		
		\item $G(D_6,R)=E(D_6,R)$.
		
		\item $G(\Phi,R,I^2)=E(\Phi,R,I^2)G(D_6,R,I^2)$.
		
		\item $z_{\alpha_1}(\xi,\zeta)\in E(\Phi,\Delta,R,I)$.
	\end{enumerate}
}
\begin{proof}
	\begin{enumerate}
		\item That holds for the arbitrary ring and the arbitrary root system of rank at least 2. See \cite{TitsCongruence}.
		
		\item That holds for the Chevalley groups of polynomial rings over Dedekind domains of arithmetic type for arbitrary root system of rank at least 2. See \cite{StavPolynom}.
		
		\item
		This is exactly the statement that the map at the {\it relative} $K_1$ groups
		$$
		K_1(D_6,R,I^2)=G(D_6,R,I^2)/E(D_6,R,I^2)\to K_1(E_7,R,I^2)=G(E_7,R,I^2)/E(E_7,R,I^2)
		$$
		is surjective. Corollary 1.7(b) of \cite{Stein71} claims that map is surjective whenever the map on corresponding usual $K_1$ groups over the ring $\tilde{R}$ is surjective, where $\tilde{R}$ is the ring such that the following diagram is a pullback 
		$$
		\xymatrix{
			\tilde{R}\ar[r]\ar[d] & R\ar[d]^{\pr}\\
			R\ar[r]^{\pr\quad} & R/I^2
		}\tp
		$$ 
		In our case, it is easy to see that the ring $\tilde{R}$ is again finitely generated, and its Krull dimension is equal to $3$, so we can apply Corollary \ref{K1iso} again.	
		
		\item 
		The reduction homomorphism
		$$
		\map{\rho_I}{E(\Delta,\Phi,R,I)}{E\left(\Phi,\Delta,R/I^2,I/I^2\right)}
		$$
		is surjective. Hence by Lemma \ref{moduloxi^2} $z_{\alpha_1}(\xi,\zeta)\in E(\Phi,\Delta,R,I)G(\Phi,R,I^2)$.
		
		Moreover, using the previous items we get
		\begin{align*}
		G(\Phi,R,I^2)=E(\Phi,R,I^2)G(D_6,R,I^2)\le E(\Phi,R,I^2)G(D_6,R)=E(\Phi,R,I^2)E(D_6,R)\le\\\le E(\Phi,\Delta,R,I)\tp
		\end{align*}
		
		Therefore, 	$z_{\alpha_1}(\xi,\zeta)\in E(\Phi,\Delta,R,I)$.
	\end{enumerate}
\end{proof}

Now, let $R$ and $I$ be arbitrary with $2\in R^*$. To prove Theorem \ref{EE} we must show that $z_{\alpha}(\xi,\zeta)\in E(\Phi,\Delta,R,I)$ for all $\alpha\in\Phi$, $\xi\in I$, and $\zeta\in R$. Firstly, it is enough to consider the case $\alpha\in \Phi\sm\Delta$, since otherwise $z_{\alpha}(\xi,\zeta)\in E(\Delta,R)$. Moreover, it is easy to see that the Weyl group  $W(\Delta)$ acts transitively on the set $\Phi\sm\Delta$, therefore, it is enough to consider the case $\alpha=\alpha_1$. Secondly, it is enough to consider the "universal situation", where $R=\Z[{1\over 2}][\xi,\zeta]$ and $I=(\xi)\unlhd R$. Indeed, if the theorem is proved for the universal case, then the required result for arbitrary $R$ and $I$ can be obtained by applying evaluation homomorphism from  $\Z[{1\over 2}][\xi,\zeta]$ to $R$. The theorem now follows from the fourth item of Lemma \ref{universal}.

\section{Acknowledgment}

I am very grateful to my advisor Nikolai Vavilov for the careful reading of the present paper and making important remarks and corrections.

\bibliographystyle{utf8gost71s}
\bibliography{english}

\begin{thebibliography}{10}
\def\selectlanguageifdefined#1{
\expandafter\ifx\csname date#1\endcsname\relax
\else\language\csname l@#1\endcsname\fi}
\ifx\undefined\url\def\url#1{{\small #1}}\else\fi
\ifx\undefined\BibUrl\def\BibUrl#1{\url{#1}}\else\fi
\ifx\undefined\BibAnnote\long\def\BibAnnote#1{}\else\fi
\ifx\undefined\BibEmph\def\BibEmph#1{\emph{#1}}\else\fi

\bibitem{BakPetrTang}
\selectlanguageifdefined{english}
\BibEmph{Bak~A., Petrov~V.~A., Tang~G.} Stability for quadratic ${K}_1$~//
  \BibEmph{K-Theory}. ---
\newblock 2003. ---
\newblock Vol.~30, no.~1. ---
\newblock p.~1--11.

\bibitem{BakTang}
\selectlanguageifdefined{english}
\BibEmph{Bak~A., Tang~G.} Stability for hermitian ${K}_1$~// \BibEmph{J. Pure
  and Appl. Algebra}. ---
\newblock 2000. ---
\newblock Vol. 150. ---
\newblock p.~107--121.

\bibitem{Bass64}
\selectlanguageifdefined{english}
\BibEmph{Bass~H.} {K}-theory and stable algebra~// \BibEmph{Publ. Math. Inst.
  Hautes \'Etudes Sci}. ---
\newblock 1964. ---
\newblock Vol.~22. ---
\newblock p.~5--60.

\bibitem{BassBook}
\selectlanguageifdefined{english}
\BibEmph{Bass~H.} Algebraic {K}-theory. ---
\newblock New York: W. A. Benjamin, Inc., 1986.

\bibitem{Bourbaki7-9}
\selectlanguageifdefined{english}
\BibEmph{Bourbaki~N.} Elements of Mathematics. {L}ie Groups and {L}ie Algebras.
  {C}hapters 7-9. ---
\newblock Springer, 2004.

\bibitem{Bourbaki4-6}
\selectlanguageifdefined{english}
\BibEmph{Bourbaki~N.} Elements of Mathematics. {L}ie Groups and {L}ie Algebras.
  {C}hapters 4-6. ---
\newblock Berlin-Heidelberg-New York: Springer-Verlag, 2008. ---
\newblock P.~300.

\bibitem{ChevalleySemiSimp}
\selectlanguageifdefined{english}
\BibEmph{Chevalley~C.} Certain schemas des groupes semi-simples~//
  \BibEmph{Sem. Bourbaki}. ---
\newblock 1960--1961. ---
\newblock p.~1--16.

\bibitem{EstesOhm}
\selectlanguageifdefined{english}
\BibEmph{D.Estes, J.Ohm}. Stable range in commutative rings~// \BibEmph{J.
  Algebra}. ---
\newblock 1967. ---
\newblock Vol.~7. ---
\newblock p.~343--362.

\bibitem{MagWandarKalVas}
\selectlanguageifdefined{english}
\BibEmph{Magurn~B.~A., {Wan dar Kallen}, Vaserstein~L.~N.} Absolute stable rank
  and Witt cancellation for noncommutative rings~// \BibEmph{Invent. Math.} ---
\newblock 1988. ---
\newblock Vol.~91. ---
\newblock p.~525--542.

\bibitem{Matsumoto69}
\selectlanguageifdefined{english}
\BibEmph{Matsumoto~H.} Sur les sous-groupes arithm\'etiques des groupes
  semisimples d\'eploy\'es~// \BibEmph{Ann. Sci. \'Ec. Norm. Sup\'er. $(4)$}.
  ---
\newblock 1969. ---
\newblock Vol.~2, no.~1. ---
\newblock p.~1--62.

\bibitem{PlotkinE7}
\selectlanguageifdefined{english}
\BibEmph{Plotkin~E.~B.} On the Stability of the ${K}_1$ -Functor for
  {C}hevalley Groups of Type ${E}_7$~// \BibEmph{J. Algebra}. ---
\newblock Vol. 210.

\bibitem{PlotProc}
\selectlanguageifdefined{english}
\BibEmph{Plotkin~E.~B.} Stability theorems for ${K}_1$ -functors for Chevalley
  groups~// Proceedings of the Conference on Nonassociative Algebras and
  Related Topics, Hiroshima. ---
\newblock 1991.

\bibitem{PlotkinStability}
\selectlanguageifdefined{english}
\BibEmph{Plotkin~E.~B.} Surjective stabilization for ${K}_1$ -functor for some
  exceptional {C}hevalley groups~// \BibEmph{J. Soviet Math.} ---
\newblock 1993. ---
\newblock Vol.~64. ---
\newblock p.~751--767.

\bibitem{atlas}
\selectlanguageifdefined{english}
\BibEmph{Plotkin~E.~B., Semenov~A.~A., Vavilov~N.~A.} Visual basic
  representations: an atlas~// \BibEmph{Internat. J. Algebra Comput.} ---
\newblock 1998. ---
\newblock Vol.~8, no.~1. ---
\newblock p.~61--95.

\bibitem{StavPolynom}
\selectlanguageifdefined{english}
\BibEmph{Stavrova~A.~K.} {C}hevalley groups of polynomial rings over {D}edekind
  domains. ---
\newblock 2018. ---
\newblock arXiv: 1812.04326 [math.KT].

\bibitem{SteinChevalley}
\selectlanguageifdefined{english}
\BibEmph{Stein~M.~R.} Generators, relations, and coverings of {C}hevalley
  groups over commutative rings~// \BibEmph{Amer. J. Math.} ---
\newblock 1971. ---
\newblock Vol.~93. ---
\newblock p.~965--1004.

\bibitem{Stein71}
\selectlanguageifdefined{english}
\BibEmph{Stein~M.~R.} Relativizing Functors on Rings and Algebraic
  ${K}$-Theory~// \BibEmph{J. Algebra}. ---
\newblock 1971. ---
\newblock Vol.~19. ---
\newblock p.~140--152.

\bibitem{SteinStability}
\selectlanguageifdefined{english}
\BibEmph{Stein~M.~R.} Stability theorems for $K_1$, $K_2$ and related functors
  modeled on {C}hevalley groups~// \BibEmph{Japan J. Math}. ---
\newblock 1978. ---
\newblock Vol.~4. ---
\newblock p.~77--108.

\bibitem{SteinMats}
\selectlanguageifdefined{english}
\BibEmph{Stein~M.~R.} Matsumoto’s solution of the congruence subgroup problem
  and stability theorems in algebraic K-theory~// Proceedings of the 19th
  Meeting Algebra Section Math. Soc. Japan. ---
\newblock 1983. ---
\newblock p.~32--44.

\bibitem{Steinberg85}
\selectlanguageifdefined{english}
\BibEmph{Steinberg~R.~G.} Some consequences of the elementary relations in
  ${SL_n}$~// \BibEmph{Contemp. Math.} ---
\newblock 1985. ---
\newblock Vol.~45. ---
\newblock p.~335--350.

\bibitem{StepComput}
\selectlanguageifdefined{english}
\BibEmph{Stepanov~A.~V.} Elementary calculus in {C}hevalley groups over
  rings~// \BibEmph{J. Prime Research in Math.} ---
\newblock 2013. ---
\newblock Vol.~9. ---
\newblock p.~79--95.

\bibitem{SusTul}
\selectlanguageifdefined{english}
\BibEmph{Suslin~A.~A., Tulenbaev~M.~S.} Stabilization theorem for the {M}ilnor
  ${K}_2$-functor~// \BibEmph{J. Soviet Math.} ---
\newblock 1981. ---
\newblock Vol.~17. ---
\newblock p.~1804--1819.

\bibitem{Taddei}
\selectlanguageifdefined{english}
\BibEmph{Taddei~G.} Normalit{\'e} des groupes {\'el\'e}mentaire dans les
  groupes de {C}hevalley sur un anneau~// \BibEmph{Contemp. Math.} ---
\newblock 1986. ---
\newblock Vol.~55. ---
\newblock p.~693--710.

\bibitem{TitsCongruence}
\selectlanguageifdefined{english}
\BibEmph{Tits~J.} Systemes g\'en\'erateurs de groupes de congruences~//
  \BibEmph{C. R. Acad. Sci., Paris, S\'er. A}. ---
\newblock 1976. ---
\newblock Vol. 283. ---
\newblock p.~693--695.

\bibitem{VasersteinStability}
\selectlanguageifdefined{english}
\BibEmph{Vaserstein~L.~N.} On the stabilization of the general linear group
  over a ring~// \BibEmph{Sb. Math.} ---
\newblock 1969. ---
\newblock Vol.~8, no.~3. ---
\newblock p.~383--400.

\bibitem{VasersteinUniStab}
\selectlanguageifdefined{english}
\BibEmph{Vaserstein~L.~N.} Stabilization of unitary and orthogonal groups over
  a ring with involution~// \BibEmph{Sb. Math.} ---
\newblock 1970. ---
\newblock Vol.~10, no.~3. ---
\newblock p.~307--326.

\bibitem{VasersteinChevalley}
\selectlanguageifdefined{english}
\BibEmph{Vaserstein~L.~N.} On normal subgroups of {C}hevalley groups over
  commutative rings~// \BibEmph{Tohoku Math. J.} ---
\newblock 1986. ---
\newblock Vol.~38. ---
\newblock p.~219--230.

\bibitem{VasSusSerre}
\selectlanguageifdefined{english}
\BibEmph{Vaserstein~L.~N., Suslin~A.~A.} {S}erre's problem on projective
  modules over polynomial rings, and algebraic {K}-theory~// \BibEmph{Math.
  USSR. Izvestija}. ---
\newblock 1976. ---
\newblock Vol.~10, no.~5. ---
\newblock p.~937--1001.

\bibitem{VavThirdlook}
\selectlanguageifdefined{english}
\BibEmph{Vavilov~N.~A.} A third look at weight diagrams~// \BibEmph{Rend.
  Semin. Mat. Univ. Padova}. ---
\newblock 2000. ---
\newblock Vol. 104. ---
\newblock p.~201--250.

\bibitem{VavSigns}
\selectlanguageifdefined{english}
\BibEmph{Vavilov~N.~A.} Can one see the signs of structure constants?~//
  \BibEmph{St.Petersburg Math. J.} ---
\newblock 2008. ---
\newblock Vol.~19. ---
\newblock p.~519--543.

\end{thebibliography}

\end{document}